\documentclass[11pt]{article}
\usepackage{amsmath,amsthm,amssymb}
\usepackage{fancyhdr}
\usepackage[british]{babel}
\usepackage{geometry,mathtools}
\usepackage{enumitem}
\usepackage{algpseudocode}
\usepackage{dsfont}
\usepackage{centernot}
\usepackage{xstring}
\usepackage{colortbl}
\usepackage[]{algorithm}
\usepackage{graphicx}
\usepackage[font={footnotesize}]{caption}
\usepackage[usenames,dvipsnames,table]{xcolor}
\usepackage{pstricks,tikz}
\usepackage[h]{esvect}
\usepackage[
		bookmarksopen=true,
		bookmarksopenlevel=1,
		colorlinks=true,
		linkcolor=darkblue,
        linktoc=page,
		citecolor=darkblue,
]{hyperref}
\usepackage{bbm}

\definecolor{darkblue}{rgb}{0,0,0.5}

\newdimen\margin
\def\textno#1&#2\par{
   \margin=\hsize
   \advance\margin by -4\parindent
          \setbox1=\hbox{\sl#1}
   \ifdim\wd1 < \margin
      $$\box1\eqno#2$$
   \else
      \bigbreak
      \hbox to \hsize{\indent$\vcenter{\advance\hsize by -3\parindent
      \it\noindent#1}\hfil#2$}
      \bigbreak
   \fi}

\newtheorem{theorem}[algorithm]{Theorem}

\theoremstyle{definition}

\newtheorem{conj}[algorithm]{Conjecture}

\def\proof{\removelastskip\penalty55\medskip\noindent\begin{stepenv}\end{stepenv}{\bf Proof. }} 

\def\noproof{{\unskip\nobreak\hfill\penalty50\hskip2em\hbox{}\nobreak\hfill%
       $\square$\parfillskip=0pt\finalhyphendemerits=0\par}\goodbreak}
\def\endproof{\noproof\bigskip}

\newcounter{stepenv}
\newenvironment{stepenv}[1][]{\refstepcounter{stepenv}}{}

\newcounter{step}[stepenv]

\newcounter{substep}[step]
\renewcommand{\thesubstep}{\thestep.\arabic{substep}}

\newcounter{claim}[stepenv]

\newcommand{\cO}{\mathcal{O}}

\newcommand{\prob}[1]{\mathrm{\mathbb{P}}\left[#1\right]}
\newcommand{\cprob}[2]{\prob{#1 \;\middle|\; #2}}
\newcommand{\expn}[1]{\mathrm{\mathbb{E}}\left[#1\right]}

\def\sm{\setminus}

\def\In{\subset}

\def\COMMENT#1{}
\def\TASK#1{}
\let\TASK=\footnote             

\begin{document}

\title{A note on dense bipartite induced subgraphs}

\author{Stefan Glock \thanks{Institute for Theoretical Studies, ETH, 8092 Z\"urich, Switzerland.
Email: \href{mailto:dr.stefan.glock@gmail.com}{\nolinkurl{dr.stefan.glock@gmail.com}}.
Research supported by Dr. Max R\"ossler, the Walter Haefner Foundation and the ETH Z\"urich Foundation.}
}

\date{}

\maketitle

\begin{abstract} 
This exposition contains a short and streamlined proof of the recent result of Kwan, Letzter, Sudakov and Tran that every triangle-free graph with minimum degree $d$ contains an induced bipartite subgraph with average degree $\Omega(\ln d/\ln\ln d)$.
\end{abstract}

\section{Introduction}

The problem of finding induced subgraphs with many edges in triangle-free graphs was initiated by Erd\H{o}s, Faudree, Pach and Spencer~\cite{EFPS:88} around 30 years ago.
Recently, Esperet, Kang and Thomass\'e posed the following interesting conjecture, where instead of finding a subgraph with many edges, one aims to maximize the average degree.

\begin{conj}[\cite{EKT:19}]
	Any triangle-free graph with minimum degree at least $d$ contains an induced bipartite subgraph with average degree $\Omega(\ln d)$.
\end{conj}

They showed that the result would be optimal up to the implied constant, by considering an appropriate binomial random graph.
Moreover, if a triangle-free graph is $d$-regular, it has chromatic number $\cO(d/\ln d)$ by Johansson's theorem, and a simple averaging argument shows that two color classes in such a coloring will induce the desired subgraph. 
However, the general problem seems to be harder. Esperet, Kang and Thomass\'e remarked that it seems to be difficult to find even an induced bipartite subgraph with minimum degree~$3$.
Kwan, Letzter, Sudakov and Tran~\cite{KLST:20} made significant progress on this problem, by showing that one can guarantee an induced bipartite subgraph with average degree~$\Omega(\ln d/\ln\ln d)$.
Their proof is already relatively short and we follow the same outline here. In essence, we omit the necessity for a more involved `cleaning procedure' by exploiting more randomness in the setup. We provide this short proof for the convenience of other researchers.

\section{Theorem and proof}

\begin{theorem}[\cite{KLST:20}]
	Any triangle-free graph with minimum degree at least $d$ contains an induced bipartite subgraph with average degree $\Omega(\ln d/\ln\ln d)$.
\end{theorem}

\proof
Define $\ell:=\lfloor\ln d/\ln\ln d \rfloor$. 
Let $G$ be a triangle-free graph on $n$ vertices with minimum degree at least~$d$, where we assume that $d$ is sufficiently large.

The first key insight is that we can reduce to the case that $G$ is $d$-degenerate. Indeed, if some proper induced subgraph $G'$ of $G$ would have minimum degree at least~$d$, then we could simply continue with $G'$ instead of $G$. Thus, we may assume that every proper induced subgraph of $G$ has a vertex with degree less than~$d$, which implies that $G$ is $d$-degenerate.

Henceforth, fix an ordering of the vertices from left to right such that every vertex has at most $d$ `left-neighbors'.
Moreover, for every vertex $v$, fix some arbitrary set $N_v\In N_G(v)$ of size exactly~$d$. 

We produce the desired subgraph as follows: Let $X$ be a random set where every vertex is included independently with probability $p:=1/d$, and let $I$ consist of the vertices in $X$ which do not have a left-neighbor in~$X$. Clearly, $I$ is independent. This set will form one part of our induced bipartite subgraph. For the other part, we need vertices that have many neighbors in~$I$. Let $Y$ be the set of vertices $y$ for which $|N_y\cap X| = \ell$, and let $Y'$ be the set of vertices in $Y$ which have at least $\ell/10$ neighbors in~$I$. Our aim is to guarantee in $Y'$ an independent set that has roughly the same size as~$I$. To this end, we show that in expectation, $Y'$ is sufficiently large, but induces very few edges.

Clearly, $\expn{|X|}=pn$.
For every vertex $y$, the random variable $|N_y\cap X|$ has binomial distribution with parameters $d,p$, hence $$\prob{y\in Y}=q:=\binom{d}{\ell}p^{\ell}(1-p)^{d-\ell}.$$
Moreover, for every edge $xy$, since $G$ is triangle-free, $N_x$ and $N_y$ are disjoint and hence the events $x\in Y$ and $y\in Y$ are independent. Thus, $\expn{e(Y)}=q^2e(G) \le q^2 nd$. 

The following are the crucial claims.
First,
\begin{align}
q\ge 0.35p. \label{crucial condition binomial}
\end{align}
This is easy to check and the very reason for our choice of~$\ell$. Indeed, we have
\begin{align*}
	q&=\binom{d}{\ell}p^{\ell}(1-p)^{d-\ell} \ge \left(\frac{d}{\ell}\right)^{\ell}p^\ell \left(1-p\right)^d \ge 0.35 \ell^{-\ell} \ge 0.35 p,
\end{align*}
where the last inequality holds because $\ell^\ell \le d$ by definition of~$\ell$.

Secondly, for every vertex $y$, we have 
\begin{align}
	\cprob{y\in Y'}{y\in Y} > 1/5. \label{crucial condition}
\end{align}

Assume that this is true.
We immediately obtain that $\expn{|Y'|}> qn/5$. Combined with our previous observations and linearity of expectation, 
$$\expn{|Y'|-e(Y)/10qd - q|X|/10p} > qn/5 - qn/10 - qn/10 = 0.$$
By the probabilistic method, there exists a choice for the set $X$ such that $|Y'|-e(Y)/10qd - q|X|/10p >0$. In particular, we deduce that $|Y'|>0$, that $e(Y')\le e(Y)\le 10qd |Y'|$ and that $|Y'|\ge q|X|/10p$.
Crucially, the average degree of $G[Y']$ is implied to be at most $20qd$. By Tur\'an's theorem, we can find an independent set $J$ inside $Y'$ of size $|J|\ge |Y'|/(20qd+1) >0$.
We claim that $G[I,J]$ is the desired subgraph. Since every vertex in $J$ has at least $\ell/10$ neighbors in $I$ by definition of $Y'$, we only have to check that $|I|\le |X|\le 10p|Y'|/q \le (200+10p/q)|J| \le 230|J|$, where the last inequality holds by~\eqref{crucial condition binomial}. Evidently, $G[I,J]$ has average degree at least~$\ell/2310$.

It remains to prove~\eqref{crucial condition}. Fix any vertex $y$. Note that the event $y\in Y$ is determined by the random choices made for $N_y$. Fix any such outcome for which $y\in Y$ holds, that is, let $L\In N_y$ have size $\ell$ and condition on the event $X\cap N_y=L$. 
Now, we reveal the remaining random choices. For a vertex $x\in L$, we have $x\in I$ if and only if no left-neighbor of $x$ is in $X$. Since $N_y$ does not contain any neighbor of $x$ by the triangle-freeness of $G$, the conditioning does not influence this event. Crucially, $x$ has at most $d$ left-neighbors, hence the probability of $x\in I$ is at least $(1-p)^d\ge 0.35$. Thus, the expected number of vertices in $L$ which are not in $I$ is at most $0.65\ell$. By Markov's inequality, the probability that $|L\sm I|\ge 0.9\ell$ is at most $\frac{0.65\ell}{0.9\ell}< 0.8$. Therefore, with probability $>1/5$, we have  $|L\cap I|\ge \ell/10$ and hence $y\in Y'$. Since this holds for any choice of $L$, the claim follows.
\endproof

\section*{Acknowledgement} Thanks to Benny Sudakov for helpful discussions leading to this note.

\providecommand{\bysame}{\leavevmode\hbox to3em{\hrulefill}\thinspace}
\providecommand{\MR}{\relax\ifhmode\unskip\space\fi MR }
\providecommand{\MRhref}[2]{%
	\href{http://www.ams.org/mathscinet-getitem?mr=#1}{#2}
}
\providecommand{\href}[2]{#2}

\end{document}